\documentclass{amsart}
\usepackage{amssymb}
\newcommand{\nor}[1]{|\!|#1|\!|}
\newcommand{\suc}{{\operatorname{\mathsf    {succ}}}}
\newcommand{\stem}{{\operatorname{\mathsf   {stem}}}}
\newcommand{\Pg}{{\mathbb P}}
\newcommand{\forces}{\Vdash}
\newcommand{\ZFCa}{{\operatorname{\mathsf {ZFC}}}}
\newcommand{\CH}{\operatorname{\mathsf {CH}}}
\newcommand{\reals}{{\mathbb R}}
\newcommand{\rationals}{{\mathbb Q}}
\newcommand{\rest}{{\mathord{\restriction}}}

\newcommand{\unif}{\operatorname{\mathsf  {non}}}

\newcommand{\cf}{{\operatorname{\mathsf   {cf}}}}
\newcommand{\dom}{{\operatorname{\mathsf {dom}}}}
\newcommand{\ran}{{\operatorname{\mathsf    {range}}}}
\newcommand{\N}{{\mathcal N}}
\newcommand{\M}{{\mathcal M}}
\newcommand{\V}{{\mathbf V}}
\newcommand{\<}{\langle}
\renewcommand{\>}{\rangle}
\newcommand{\QED}{\hspace{0.1in} \square \vspace{0.1in}}
\newcommand{\thinks}{\models}
\newcommand{\Proof}{{\sc Proof} \hspace{0.2in}}
\newcommand{\lft}[2]{\mathopen\ifcase#1{}\oo\or
                        \big#2\or\Big#2\else\oo\fi} 
\newcommand{\rgt}[2]{\mathclose\ifcase#1{}\oo\or
                        \big#2\or\Big#2\else\oo\fi} 

\newcommand{\SN}{{\mathcal  {SN}}}

\newcommand{\UNIF}{{\mathsf   {NON}}}

\theoremstyle{plain}
\newtheorem{theorem}{Theorem}
\theoremstyle{plain}
\newtheorem{lemma}[theorem]{Lemma}

\newtheorem{definition}[theorem]{Definition}

\begin{document}
\title{Continuous images of sets of reals}
\author{Tomek Bartoszynski}
\address{Department of Mathematics and Computer Science\\
Boise State University\\
Boise, Idaho 83725 U.S.A.}
\thanks{The first author was  partially supported by 
NSF grant DMS 9971282} 
\email{tomek@math.idbsu.edu, http://math.idbsu.edu/\char 126 tomek}
\author{Saharon Shelah}
\thanks{The second author was partially supported by Israel Science
   Foundation. Publication 722}
\address{Department of Mathematics\\
Hebrew University\\
Jerusalem, Israel}
\email{shelah@sunrise.huji.ac.il, http://math.rutgers.edu/\char 126 shelah/}
\begin{abstract}
We will show that, consistently, every uncountable set can be
continuously mapped onto a non measure zero set, while there exists an 
uncountable set whose all continuous images into a Polish space are meager.
\end{abstract}
\maketitle

\section{Introduction}

Let $ {\mathcal J} $ be a $ \sigma $-ideal of subsets
of a Polish space $Y$. Assume also that ${\mathcal J} $ contains
singletons and has Borel basis.
Let $\unif({\mathcal J})=\min\{|X| : X \subseteq Y\ \&\ X \not \in {\mathcal J} \}.$

In this paper we are concerned with the family 
$$\UNIF({\mathcal J})=\{X
\subseteq \reals : \text{for every continuous mapping } F:X
\longrightarrow Y,\ F"(X) \in {\mathcal J}\}.$$
Note that
$\UNIF({\mathcal J})$ contains all countable sets. 
Moreover, $\UNIF({\mathcal J})$ is closed under countable unions but
need not be downward closed, thus it may not an ideal.
However, $\UNIF({\mathcal J})$ is contained in the $ \sigma $-ideal
$$\UNIF^\star({\mathcal J})=\{X
\subseteq \reals : \text{for every continuous mapping } F:\reals
\longrightarrow Y,\ F"(X) \in {\mathcal J}\}.$$
It is also not hard to see that  $\UNIF^\star({\mathcal J})$ consists
of those sets whose uniformly continuous images are in $ {\mathcal J}
$.

Let $\N$ be the ideal of measure zero subsets of $2^\omega $ with
respect to the standard product measure $\mu$, and let $\M$ be the
ideal of meager subsets of $2^\omega $ (or other Polish space $Y$).

We will show that
\begin{itemize}
\item $\ZFCa \vdash \UNIF(\M)$  contains an uncountable set.
\item It is consistent that  $\UNIF^\star(\N)=\UNIF(\N)=[\reals]^{\leq
  \mathbf\aleph_0}$.
\item It is consistent that:  $\UNIF^\star({\mathcal J})=[\reals]^{\leq
  \mathbf\aleph_0} \iff \unif({\mathcal J})<2^{\boldsymbol\aleph_0} $.
\end{itemize}
Observe that if $\UNIF^\star({\mathcal J})=[\reals]^{\leq
  \mathbf\aleph_0}$ then $\unif({\mathcal
  J})=\boldsymbol\aleph_1 $. On the other hand for all $ \sigma $-ideals
${\mathcal J} $ considered in this paper, continuum hypothesis implies 
that $\UNIF({\mathcal J})$ contains uncountable sets. 

Finally notice that one can show in $\ZFCa$ that there exists an
uncountable universal measure zero set (see \cite{Mil84Spe}), i.e. a 
set whose all homeomorphic (or even Borel isomorphic)
images are all of measure zero. Therefore one cannot generalize the
consistency results mentioned above by replacing the word
``continuous'' by ``homeomorphic'' in the definition of $\UNIF(\N)$.
 
\section{Category}
In this section we show that $\UNIF(\M)$ contains an uncountable
set. This was proved in \cite{JMSS}, the proof presented here gives a
slightly stronger result.

For $f,g \in \omega^\omega $ let $f \leq^\star g$ mean that $f(n) \leq 
g(n)$ for all but finitely many $n$.
Let 
$${\mathfrak b}=\min\{|F|: F \subseteq \omega^\omega \ \&\ \forall h
\in \omega^\omega \ \exists f \in F \ f \not \leq^\star h\}.$$

\begin{theorem}[\cite{JMSS}]\label{1}
There exists a set $X \subseteq \reals$ of size ${\mathfrak b}$ 
such that
\begin{enumerate}
\item 
every continuous image of $X$ into $\omega^\omega $ is bounded,
\item every continuous image of $X$ into a Polish space is meager,
\item if ${\mathfrak b} \leq \unif(\N)$ then every continuous image of
  $X$ into $\reals$ has measure zero.
\end{enumerate}
\end{theorem}
\Proof
Let $Z \subseteq (\omega+1)^\omega $ consist of functions $f$ such
that 
\begin{enumerate}
\item $\forall n \ f(n) \leq f(n+1)$,
\item $\forall n\ \lft1(f(n)<\omega \rightarrow f(n)<f(n+1)\rgt1)$.
\end{enumerate}
Note that  $Z$ is a compact subset of $(\omega+1)^\omega $ thus it is
homeomorphic to $2^\omega $.
For an increasing sequence $s \in (\omega+1)^{<\omega}$ let $q_{s} \in (\omega+1)^\omega $ be
defined as
$$q_{s}(k)=\left\{
  \begin{array}{ll}
s(k)& \text{if } k <|s|\\
\omega & \text{otherwise}
  \end{array}\right. \text{ for } k \in \omega . $$
Note that the set  $\rationals=\{q_s: s \in \omega^{<\omega} \}$
is dense in $Z$.
Put $X'=\{f_\alpha : \alpha < {\mathfrak b}\}$ such that 
\begin{enumerate}
\item $ \forall \alpha \ f_\alpha \in Z$,
\item $f_\alpha \leq^\star f_\beta $ for $\alpha<\beta $,
\item $\forall f \in \omega^\omega \ \exists \alpha \ f_\alpha \not
  \leq^\star f$. 
\end{enumerate}

Let $X = X' \cup \rationals$. We will  show that $X$ is the set we are 
looking for.

(1) Suppose that
$F:X \longrightarrow \omega^\omega $ is continuous. We only need to
assume that $F$ is continuous on $\rationals$. Without loss of
generality we can assume for every $x \in X$, $F(x) \in \omega^\omega
$ is an increasing function.
\begin{lemma}
  There exists a function $g \in \omega^\omega $ such that for every
  $x \in X$ and $n \in \omega $,
$$F(x)(n)\leq g(n) \text{ if } x(n)>g(n).$$
\end{lemma}
\Proof
Fix $n \in \omega $ and for each $s \in (\omega+1)^n$ let $I_s
\subseteq (\omega+1)^\omega $ be a
basic open set containing $q_s$ such that for $x \in \dom(F)\cap I_s$,
$F(x)(n)=F(q_s)(n)$.
For every $s$ the set $I_s \rest n=\{x\rest n: x \in I_s\}$ is open
$(\omega+1)^n$ and the family $\{I_s\rest n: s\in (\omega+1)^n\}$ is a 
cover of $(\omega+1)^n$. By compactness there are sequences $s_1,
\dots, s_k$ such that 
$(\omega+1)^n=I_{s_1}\rest n \cup \dots\cup I_{s_k}\rest n.$
Find $N$ so large that if $x(n)>N$ then $x \in I_{s_j}$ for some
$j\leq k$.
Define
$$g(n)=\max\{N,F(q_{s_1})(n), \dots, F(q_{s_k})(n)\}.~\QED$$

Let $g \in \omega^\omega $ be the function from the above lemma. Find
$\alpha_0$ such that $f_{\alpha_0} \not \leq^\star g$.
Let $\{u_n: n \in \omega\}$ be an increasing enumeration of $\{n:
g(n)<f_{\alpha_0}(n)\}$.
Put $h(n)=g(u_n)$ for $n \in \omega $ and note that for $
\beta>\alpha_0$ and sufficiently large $n$ we have
$$F(f_\beta)(n)\leq F(f_\beta)(u_n) < g(u_n)=h(n).$$
Since the set $\{F(f_\beta): \beta \leq \alpha_0\} \cup \{F(q_s): s
\in \omega^{<\omega}\}$ has size $< {\mathfrak b}$ we conclude that
$F"(X)$ is bounded in $\omega^\omega $.

\bigskip

(2) Suppose that $F$ is a continuous mapping from $X$ into a Polish
space $Y$ with metric $\rho$. Observe that $F$ is not onto and fix a countable dense set
$\{q_n: n \in \omega \}$ disjoint with $F"(X)$. 
For $x \in X$ let $f_x \in \omega^\omega $ be defined as
$$f_x(n)=\min\left\{k: \rho(f(x),d_n)>\frac{1}{k}\right\}.$$
In particular, 
$$f(x) \not \in B\left(d_n, \frac{1}{f_x(n)}\right)=\left\{z: \rho(d_n,z)<\frac{1}{f_x(n)}\right\}.$$
Note that the mapping $x \mapsto f_x$ is continuous and find a
function $h \in \omega^\omega $ such that 
$f_x \leq^\star h$ for  $x \in X$.
Put
$$G=\bigcap_m\bigcup_{n>m} B\left(d_n, \frac{1}{h(n)}\right)$$
and
note that $G$ is a comeager set disjoint from $F"(X)$.

\bigskip

(3) Let $\rationals \subseteq U$ be an open set.  Define $g \in
\omega^\omega $ as
$$g(0)=\min\{k: \forall x\ x(0)>k \rightarrow x \in U\}$$
and for $n>0$
$$g(n)=\min\left\{k: \forall x \ \lft2(\lft1(\forall j<n \ x(j)<g(j) \
  \& \ x(n)>k\rgt1) 
\rightarrow x \in U\rgt2)\right\} .$$
Let $\alpha_0$ be such that $f_{\alpha_0}\not \leq^\star g$. It
follows that $f_\beta \in U$ for $\beta>\alpha_0$.

Suppose that $F:X \longrightarrow \reals$ is continuous (on
$\rationals$).  Let $\{q_n: n \in \omega\}$ be enumeration of
$\rationals$.  Let $I_n^k \ni q_n$ be a basic open set such that
$F"(I^k_n)$ has diameter $<2^{-n-k}$. Put $H=\bigcap_k\bigcup_{n}
I_n^k$.  It is clear that $F"(H)$ has measure zero.  Fix $\alpha_0$
such that for all $\beta>\alpha_0$, $f_\beta \in H$. It follows that
for $\beta>\alpha_0$, $F(f_\beta)$ belongs to a measure zero set
$F"(H)$. By the assumption, the remainder of the set $F"(X)$ has size
$<\unif(\N)$ which finishes the proof.~$\QED$

The set $X$ constructed above is not hereditary, and for example $X
\setminus \rationals$ can be continuously mapped onto an unbounded
family.
A hereditary set having property (1) of theorem \ref{1} cannot be constructed
in $\ZFCa$. Miller showed in \cite{Mil79Len} that it is consistent that every
uncountable set has a subset that can be mapped onto an unbounded
family. This holds in a model where there are no $ \sigma $-sets,
i.e. every uncountable set has a $G_\delta $ subset which is not
$F_\sigma $. 

It is open whether a hereditary set having property (2) of theorem \ref{1} can
be constructed in $\ZFCa$.

\section{Making $\UNIF({\mathcal J})$ small.}

We will start with the following:
\begin{definition}
  A set $X \subseteq 2^\omega $ has strong measure zero if for every
  function $g \in \omega^\omega $ there exists a function $f \in
  (\omega^{<\omega})^\omega $ such that $f(n)\in 2^{g(n)}$ for every
  $n$ and
$$\forall x \in X\ \exists^\infty n \ x \rest g(n)=f(n).$$
Let $\SN$ denote the collection of all strong measure zero sets.

If the above property fails for some $g$ then we say that $g$ witnesses 
that $X\not\in \SN$.
\end{definition}

For $g \in \omega^\omega $ we will define a forcing notion ${\mathbb
  P}_g$ such that:
\begin{enumerate}
\item ${\mathbb P}_g$ is proper,
\item there exists a family $\{F_n: n \in \omega\}\in \V^{{\mathbb P}_g}$ such that 
  \begin{enumerate}
  \item $\forall n \ F_n:2^\omega \longrightarrow 2^\omega $ is
    continuous,
  \item if  $X \subseteq 2^\omega $, $X \in \V$ and $g$ witnesses that
    $X\not\in \SN$ then $$\V^{\Pg_g} \thinks \bigcup_n \
    F_n"(X)=2^\omega.$$
  \end{enumerate}
\end{enumerate}
Let ${\mathbb L}$ be the Laver forcing and suppose that $g$ is a Laver
real over $\V$. It is well known that:

\begin{lemma}\cite{Lav76Con}
If $X \subseteq 2^\omega $, $X \in \V$ is uncountable then
  $\V[g] \thinks $''$g$ witnesses that $X$ does not have
  strong measure zero.''
\end{lemma}

\begin{theorem}\label{main}
  It is consistent with $\ZFCa$ that for every $\sigma$-ideal ${\mathcal J}$ 
  $$\UNIF({\mathcal J} )=[\reals]^{\leq
    {\mathbf\aleph}_0}\iff \unif({\mathcal J})<2^{\boldsymbol\aleph_0}.$$
\end{theorem}
\Proof
Let $\<{\mathcal P}_\alpha, \dot{{\mathcal
    Q}}_\alpha:\alpha<\omega_2\>$ be a countable support iteration
such that $\forces_\alpha \dot{{\mathcal Q}}_\alpha \simeq {\mathbb
  L}\star {\mathbb P}_{\dot{g}}$ for $\alpha < \omega_2$.
Suppose that ${\mathcal J} $ is a $\sigma $-ideal and $\V^{{\mathcal
    P}_{\omega_2}} \thinks \unif({\mathcal J})=\boldsymbol\aleph_1 $. 
It follows that for some $\alpha<\omega_2$, $\V^{{\mathcal
    P}_{\omega_2}} \thinks \V^{{\mathcal P}_{\alpha}}\cap 2^\omega \not \in {\mathcal J}
$. 

Suppose that $X \subseteq \V^{{\mathcal
    P}_{\omega_2}} \cap 2^\omega $ is uncountable. Let $\beta >\alpha$
 be such that $X \cap \V^{{\mathcal
    P}_{\beta  }}$ is uncountable. In $\V^{{\mathcal
    P}_{\beta  }\star {\mathbb L}}$ the Laver real  witnesses that 
$X\not\in \SN$ and so $\V^{{\mathcal
    P}_{\beta+1}} \thinks \bigcup_n F_n"(X) =2^\omega  $. Hence
$\V^{{\mathcal
    P}_{\omega_2}} \thinks \bigcup_n F_n"(X) \not \in {\mathcal J} $
which means that $\V^{{\mathcal
    P}_{\omega_2}} \thinks \exists n \in \omega  \ F_n"(X) \not \in {\mathcal J}$.~$\QED$

\section{Defintion of ${\mathbb P}_g$}
Let us fix the following notation. 
Suppose that $\<F_n:n \in \omega\>$ are nonempty sets. Let
$T^{max}=\bigcup_n \prod_{j=0}^{n-1} F_j$. 
For a tree $T \subseteq T^{max}$ let $T \rest n = T \cap
\prod_{j=0}^{n-1} F_j$.
For $t \in T\rest n$
let $\suc_T(t)=\{x \in F_n: t^\frown x \in T\}$  be the set of
all immediate successors of $t$ in $T$, and 
let $T_t=\left\{s \in T : s \subseteq t \hbox{ 
or } t \supseteq s\right\}$ be the subtree determined by $t$. 
Let 
$\stem(T)$ be the 
shortest $t\in T$ such that $|\suc_T(t)|>1$.

Fix a sequence $\<\varepsilon^k_j: j \leq 
k\} $ such that 
\begin{enumerate}
\item $\forall k \ 0<\varepsilon^k_0<\varepsilon^k_1 <\dots<
  \varepsilon^k_k<2^{-k}$.
\item $  \forall k \ \forall j < k \ (\varepsilon^k_{j+1})^3 > \varepsilon^k_j$.

\item $  \forall k \ \forall j < k \ \frac{\varepsilon^k_{j+1}}{2^{k^2}} > \varepsilon^k_j$.
\item $\forall k \ \forall j<k \
  \frac{\varepsilon^k_j}{\varepsilon^k_{j+1}} < \varepsilon^k_k$.
\end{enumerate}
For example $\varepsilon^k_j=2^{-k^2 4^{k-j}}$ for $j \leq k$ will
work.

Suppose that a strictly increasing function $g \in \omega^\omega $ is given.
Fix an increasing sequence $\<n_k: k\in \omega \>$ such that
$n_0=0$ and
$$n_{k+1} \geq g\left((k+1)2^{n_k}\frac{2^{n_k 2^{n_k}}}{\varepsilon^k_0}\right) \text{ for }
k \in \omega .$$
For the choice of $\varepsilon^k_j $ above
$n_{k+1}=g\left(2^{1999^{n_k}}\right)$ will be large enough.

Let $F_k=\{f:\dom(f)= 2^{[n_k,n_{k+1})}, \ran(f)=\{0,1\}\}$.
For $A \subseteq F_k$ let 
$$\nor{A}=\max\left\{\ell: \dfrac{|A|}{|F_k|} \geq
  \varepsilon^k_\ell\right\}.$$
Consider the tree
$$T^{max}=\bigcup_k \prod_{j=0}^k F_j.$$
Let $\Pg_g $ be the forcing notion which consists of perfect subtrees $T
\subseteq T^{max}$ such that 
$$\lim_{k \rightarrow \infty} \min\{\nor{\suc_T(s)}: s \in T\rest
k\}=\infty .$$
For $T, S \in \Pg_g $ and $ n \in \omega $ define $T \geq S$ if 
$T \subseteq S$ and $T \geq_n S$ if
$T \geq S$ and
$$\forall s \in S \ \lft1(\nor{\suc_S(s)}\leq n \rightarrow \suc_S(s)=\suc_T(s)\rgt1).$$
It is easy to see that $\Pg_g$ satisfies Axiom A, thus it is proper.

Suppose that $G \subseteq \Pg_g $ is a generic filter over $\V$. Let
$\bigcap G=\<f_0,f_1,f_2, \dots\> \in \prod_k F_k$. 
Define
$F_G: 2^\omega \longrightarrow 2^\omega $ as
$$F_G(x)(k)=f_k\lft1(x\rest[n_k,n_{k+1})\rgt1) \text{ for } x\in 2^\omega ,\ k
\in \omega .$$

First we show that $\Pg_g$ is $\omega^\omega $-bounding. The arguments 
below are rather standard, we reconstruct them here for completeness
but the reader familiar with \cite{RoSh470} will see that they are a
part of a much more general scheme.

\begin{lemma}\label{5}
  Suppose that $I \subseteq \V$ is a countable set, $n \in \omega $ and $T
  \forces_{\Pg_g} \dot{a} \in I$.
There exists $S \geq_n T$ and $k \in \omega $ such that for every $t
\in S \rest k$ there exists $a_t \in I$ such that $S_t
\forces_{\Pg_g} \dot{a}=a_t$.
\end{lemma}
\Proof
Let $S \subseteq T$ be the set of all $t \in  T$ such that
$T_t$ satisfies the lemma. In other words
$$S = \{ t  \in T: \exists k_{t} \in \omega \
\exists T' \geq_{n} T_t \
\forall s \in T'\rest k_t  \ \exists a_{s} \in I
\ T'_s \forces_{\Pg_g}  \dot{a} =a_{s}\} .$$
We want to show that $\stem(T)\in S$.
Notice that if $s \not \in S$ then
$$\nor{\suc_S(s)} \leq \varepsilon^{|s|}_n.$$

Suppose that $\stem(T)\not\in S$ and by induction on levels
build a tree $\bar{S} \geq_{n} T$ such that for
$s \in \bar{S}$,
$$\suc_{\bar{S}}(s) = \left\{ \begin{array}{ll} \suc_{T}(s) & \hbox{if }
\nor{\suc_{T}(s)}\leq n \\
\suc_{T}(s) \setminus \suc_{S}(s) & \hbox{otherwise}
\end{array} \right. .$$
Clearly $\bar{S} \in {\Pg_g} $ since $\nor{\suc_{\bar{S}}(s)}\geq 
\nor{\suc_{\bar{T}}(s)}-1$ for $s$ containing $\stem(T)$.
That is a contradiction since  $\bar{S}\cap S=\emptyset$
which is impossible.~$\QED$

In our case we have even stronger fact:
\begin{lemma}\label{6}
  Suppose that $T \forces_{\Pg_g} \dot{A} \subseteq
  2^{<\omega}$. There exists $S \geq T$ such that for all but finitely 
  many $n$, for every $t \in S
  \rest n$ there exists $A_t \subseteq 2^{n}$ such that $S_t
  \forces_{\Pg_g} \dot{A}\rest n=A_t$. 

In particular, if $T \forces_{\Pg_g} \dot{x}\in 2^\omega $then
there
 exists $S \geq T$ such that for every for all but finitely 
  many $n$, for every $t \in S
  \rest n$ there exists $s_t \in 2^{n}$ such that $S_t
  \forces_{\Pg_g} \dot{x}\rest n=s_t$. 
\end{lemma}
\Proof
It is enough to prove the first part.
By applying lemma \ref{5} we can assume that there exists an increasing
sequence $\<k_n: n \in \omega\>$ such that 
for every $t \in T
  \rest k_n$ there exists $A_t \subseteq 2^{<n}$ such that $T_t
  \forces_{\Pg_g} \dot{A}\rest n=A_t$. 

Let $n_0=|\stem(T)|$. 
Build by induction a family of trees $\{T_{n,l} : n>n_0, n \leq
l       \leq 
k_{n}\}$ such that
$$\forall s \in T_{n,l} \rest l \ \exists A_{s} \subseteq 2^n \
(T_{n,l})_s \forces_{\Pg_g}  \dot{A}\rest n=A_{s}.$$


Let $T_{n_0+1,k_{n_0}}=T$ and suppose that $T_{n,l}$ has been
  constructed. If $l=n$ let $T_{n+1,k_{n+1}}=T_{n,n}$, otherwise
  construct $T_{n,l-1}$ as follows -- by the induction hypothesis for
  $s \in T_{n,l} \rest l-1$ and every $f \in \suc_{T_{n,l}}(s)$,
there exists $A_{s^\frown f} \subseteq 2^{n}$ such that 
$$(T_{n,l})_{s^\frown f} \forces_{\Pg_g} \dot{A}\rest n=A_{s^\frown f}.$$
Fix $A$ such that $\{f: A_{s^\frown f}=A\}$ has the largest size and
put
$$\suc_{T_{n,l-1}}(s)=\left\{
  \begin{array}{ll}
\{f \in \suc_{T_{n,l}}(s): A_{s^\frown f}=A\} & \text{if } s \in
T_{n,l} \rest l-1\\
\suc_{T_{n,l}}(s)&\text{otherwise}
  \end{array}\right. $$
Finally let $S =\bigcap_n T_{n,n}$. Clearly $S$ has the required
property provided that it is a member of ${\Pg_g} $.
Note that for an element $s \in S \rest k$,
$$|\suc_S(s)| \geq \frac{|\suc_{T_{k,k}}(s)|}{2^{k^2}}.$$
By the choice of sequence $\<\varepsilon^k_l:k,l\>$, it follows that
$\nor{\suc_S(s)} \geq \nor{\suc_T(s)}-1$ if $|s|>|\stem(S)|$.
Thus $S \in {\Pg_g} $ which finishes the proof.~$\QED$

Next we show that $\Pg_g$ adds a continuous function which maps sets
that do not have strong measure zero onto sets that are not in
${\mathcal J} $.

Let $\rationals=\{x\in 2^\omega : \forall^\infty n \ x(n)=0\}$ be the
set of rationals in $2^\omega $.

\begin{theorem}
  Suppose that $g$ witnesses that $X \subseteq \V\cap 2^\omega $, $X
  \in \V$ does not have strong measure zero.
Then $\forces_{\Pg_g} F_{\dot{G}}"(X)+\rationals=2^\omega $.

In particular, 
$$V^{\Pg_g} \thinks \exists  q \in \rationals \ F_q"(X) \not \in
{\mathcal J},$$
where $F_q:2^\omega \longrightarrow 2^\omega $ is 
defined as $F_q(x)=F_{\dot{G}}(x)+q$ for $x \in 2^\omega $. 
\end{theorem}
\Proof
We start with the following:
\begin{lemma}\label{4}
  Suppose that $\frac{1}{2}>\varepsilon >0$, $I \subseteq \omega $ is 
  finite and $A \subseteq 2^{2^I}$, $\dfrac{|A|}{|2^{2^I}|}\geq
  \varepsilon $.
Let
$$Z=\left\{s \in 2^I: \exists i_s\in \{0,1\}\ \frac{|\{f\in A:
    f(s)=i_s\}|}{|2^{2^I}|}<\varepsilon^3\right\}.$$ 
Then $|Z|\leq \dfrac{1}{\varepsilon}$.
\end{lemma}
\Proof
Suppose otherwise. By passing to a subset
we can assume that $|Z|=\dfrac{1}{\varepsilon}$.
Let 
$$A'=\{f\in A: \forall s \in Z\ f(s)=i_s\}.$$
By the assumption 
$$\dfrac{|A'|}{|2^{2^I}|}\geq \varepsilon - \frac{1}{\varepsilon}\cdot 
\varepsilon^3=\varepsilon-\varepsilon^2> \frac{\varepsilon}{2}.$$
On the other hand the sets $I_s=\{f \in 2^{2^I}: f(s)=0\},\ s \in 2^I$ 
are probabilistically independent and have ``measure'' $\dfrac{1}{2}$.
It follows that
$$\dfrac{|A'|}{|2^{2^I}|}\leq \frac{1}{2^{\frac{1}{\varepsilon}}} < \frac{\varepsilon}{2},$$
which gives a contradiction.~$\QED$

\begin{lemma}\label{7}
  Suppose that $T\forces_{\Pg_g} \dot{z} \in 2^\omega $. There
  exists a sequence $\<J_k: k \in \omega\>$ such that for every $k \in 
  \omega $
  \begin{enumerate}
  \item $J_k \subseteq 2^{[n_k,n_{k+1})}$,
  \item $|J_k| \leq \dfrac{2^{n_k 2^{n_k}}}{\varepsilon^k_0}$,
  \end{enumerate}
and if $x \in \V\cap 2^\omega $ and $x \rest [n_k,n_{k+1}) \not \in J_k$ 
for all but finitely many $k$ then there exists $S \geq T$ such that 
$$S \forces_{\Pg_g} F_{\dot{G}}(x) =^\star \dot{z}.$$
\end{lemma}
\Proof
Suppose that $T \forces_{{\Pg_g}} \dot{z} \in 2^\omega $.
Let $k_0=|\stem(T)|$. By lemma \ref{6}, we can
assume that
$$\forall k>k_0\ \forall t \in T \rest k\ \exists i_t \in \{0,1\} \ T_t
\forces_{\Pg_g} \dot{z}(k)=i_t.$$
For $k > k_0$ and $s \in T \rest k$ let
$$J^s_k=\{x \in 2^{[n_k, n_{k+1})}: \exists i\in \{0,1\}\ \nor{\{f\in
  \suc_T(s):f(x)=i\}}< \nor{\suc_T(s)}-1\}.$$
By lemma \ref{4}, $|J^s_k| \leq \frac{1}{\varepsilon^k_0}$.
Put 
$J_k=\bigcup_{s \in T\rest k}J^s_k$ and note that
$$|J_k| \leq \frac{1}{\varepsilon^k_0} \prod_{i=0}^{k-1}
2^{2^{n_{i+1}-n_i}} \leq \frac{2^{n_k 2^{n_k}}}{\varepsilon^k_0}.$$

Suppose that $x \rest [n_k, n_{k+1}) \not\in J_k$ for $k\geq k^\star
\geq k_0$.
Define $S \geq T$ 
$$\suc_{S}(t)=\left\{
  \begin{array}{ll}
\{f \in \suc_{T}(t): f(x\rest [n_k,n_{k+1})=i_{t}\} &
\text{if } s \in T \rest k \text{ and } k>k^\star\\
\suc_{T}(s)&\text{otherwise}
  \end{array}\right. .$$
By the choice of $x$, $\nor{\suc_S(s)} \geq \nor{\suc_T(s)}-1$ for $
s\in S$. Thus $S \in {\Pg_g} $, and 
$$S \forces_{\mathcal P} \forall k>k^\star \
F_{\dot{G}}(x)(k)=\dot{z}(k).~\QED$$

Suppose that $T \forces_{\Pg_g} \dot{z} \in 2^\omega $ and let $\{J_k: 
k \in \omega\}$ be the sequence from lemma \ref{7}. 
Let 
$$U=\{s \in 2^{<\omega}: \exists k \ |s|=n_{k+1} \ \&\ s \rest [n_k,
n_{k+1}) \in J_k\}.$$
Let $s_1, s_2, \dots $ be the list of elements of $U$ according to
increasing length. Note that by the choice of $g$, $|s_k| \geq g(k)$
for $k \in \omega $.
Since $g$ witnesses that $X \not\in \SN$ (and any bigger function
witnesses that as well) there is $x \in X$ such that 
$$\forall^\infty k \ x \rest \dom(s_k)\neq s_k.$$
Since initial parts of $s_k's$ exhaust all possibilities it follows
that for sufficiently large $k\in \omega $,
$$x \rest [n_k, n_{k+1}) \neq s_l \rest [n_k, n_{k+1}) \text{ for all
  $l$ such that $|s_l|=n_{k+1}$}.$$
In particular,
$$\forall^\infty k \ x \rest [n_k,n_{k+1}) \not\in J_k.$$
By lemma \ref{7} we conclude  that $\forces_{\Pg_g} F_{\dot{G}}(x)
=^\star \dot{z}$. Since $\dot{z}$ was arbitrary, it follows that $\forces_{\Pg_g}
F_{\dot{G}}"(X)+\rationals=2^\omega $.
As ${\mathcal J} $ is a $ \sigma $-ideal we conclude that 
$$\forces_{\Pg_g}
\exists q \in \rationals\ F_{\dot{G}}"(X)+q \not \in {\mathcal J} .~\QED$$

\section{Measure}
Theorem \ref{main} is significant only if in the constructed model
there are some interesting $  \sigma $-ideals $ {\mathcal J} $ such
that $\unif({\mathcal J})<2^{\boldsymbol\aleph_0} $.
We will show some examples of such ideals, the most important being the 
ideal of measure zero sets $\N$.
\begin{definition}
A family $\mathcal A \subseteq [\omega]^\omega$ is called a splitting family if
for every infinite 
set $B \subseteq \omega$ there exists $A \in
\mathcal A$ such that 
$$|A \cap B|=|(\omega \setminus A) \cap B|=\boldsymbol\aleph_0.$$
We say  that ${\mathcal A} $ is strongly non-splitting if for every $B 
\in [\omega]^\omega $ there exists $C \subseteq B$ which witnesses
that ${\mathcal A} $ is not splitting.
 \end{definition}
Let
$${\mathsf S} = \{X \subseteq [\omega]^\omega : X \text{ is strongly non-splitting}\}.$$
It is easy to see that ${\mathsf S}$ is a $ \sigma $-ideal.

\begin{theorem}\label{fin}
  It is consistent that for every uncountable set $X \subseteq
  2^\omega $ there exists a continuous function $F:2^\omega
  \longrightarrow 2^\omega $ such that $F"(X)$ does not have measure
  zero.
In particular, it is consistent that
$$\UNIF(\N)=\UNIF({\mathsf S})=\SN=[\reals]^{\leq{\mathbf\aleph}_0}.$$
\end{theorem}
\Proof
Let $\V^{{\mathcal P}_{\omega_2}}$ be the model constructed in the
proof of theorem \ref{main}. To show the first part it is enough to show that
$\V^{{\mathcal P}_{\omega_2}} \thinks \unif(\N)=\boldsymbol\aleph_1 $.
Since is well known that $\unif({\mathsf S})\leq \unif(\N)$,
it follows that $ \UNIF({\mathsf S})=[\reals]^{\leq{\mathbf\aleph}_0}.$
This was known to be consistent (see \cite{BarSplit}).
Finally, it is well  known that if $X \in \SN$ and $F:X
\longrightarrow 2^\omega $ is uniformly continuous then $F"(X) \in
\SN \subseteq \N$. Thus $\SN=[\reals]^{\leq{\mathbf\aleph}_0}$.

To finish the proof we have to show that $\V^{{\mathcal P}_{\omega_2}} 
\thinks \unif(\N)=\boldsymbol\aleph_1 $.
By theorem 6.3.13 of \cite{BJbook}, in order to show that it suffices to show that both
$\Pg_g$ and ${\mathbb L}$ satisfy certain condition (preservation of
$\sqsubseteq^{random}$) which is an iterable version of preservation
of outer measure.
Theorem 7.3.39 of \cite{BJbook} shows that $\mathbb L$ satisfies this
condition. Exactly the same proof works for $\Pg_g$ provided that we
show:
\begin{theorem}
  If $X \subseteq 2^\omega $, $X \in \V$ and $\V \thinks X \not
  \in \N$ then $\V^{{\mathbb P}_g} \thinks X \not
  \in \N$.
\end{theorem}
\Proof
The sketch of the proof presented here is a special case of a more
general theorem (theorem 3.3.5 of \cite{RoSh470}). 

Fix $1>\delta>0$ and a strictly increasing sequence $\<\delta_n: n \in \omega\>$ of real
numbers such that 
\begin{enumerate}
\item $\sup_n \delta_n=\delta$.
\item $\forall^\infty n \ \delta_{n+1}-\delta_n > \varepsilon^n_n$.
\end{enumerate}

Suppose that $\forces_{\Pg_g} X \not\in \N$. Without loss of
generality we can assume that $X$ is forced to have outer measure one.
Let $\dot{A}$ be a $\Pg_g$-name such that $\forces_{\Pg_g} \dot{A}
\subseteq 2^{<\omega}\ \&\ \mu([A])\geq \delta$ and suppose that
$T \forces_{\Pg_g} X \cap [\dot{A}]=\emptyset$. Let $n_0=|\stem(T)|$. 
By lemma \ref{6}, we can assume that
$$\forall n>n_0 \ \forall t \in T
  \rest n\ \exists A_t \subseteq 2^{n}\ T_t
  \forces_{\Pg_g} \dot{A}\rest n=A_t.$$

Fix $n>n_0$ and define by induction sets $\{A^n_t: t \in T\rest m, \
n_0\leq m\leq n+1\}$ such that 
\begin{enumerate}
\item $A^n_t \subseteq 2^{n+1}$ for $t \in T$,
\item $|A^n_t|\cdot 2^{-n-1} \geq \delta_m$ for $t \in T\rest m$.
\end{enumerate}

For $t \in T\rest n+1 $ let $A^n_t=A_t$. Suppose that sets $A^n_t$ are 
defined for $t \in T\rest m$, $m>n_0$.
Let $t \in T\rest m-1$ and consider the family $\{A^n_{t^\frown f}: f
\in \suc_T(t)\}$. By the induction hypothesis, $|A^n_{t^\frown f}|\cdot 2^{-n-1} \geq \delta_m$
Let
$$A^n_t=\{s \in 2^{n+1}: \nor{\{f: s \in A_{t^\frown f}} \geq
\nor{\suc_T(t)}-1\}.$$
A straightforward computation (recall Fubini theorem) shows that the requirement that we put
on the sequence $\<\delta_n: n \in \omega\>$ implies that
$|A^n_t|\cdot 2^{-n-1} \geq \delta_{m-1}$. In particular,
$A^n_{\stem(T)}\cdot 2^{-n-1} \geq \delta_{n_0}$  for all $n$.
Let $B=\{x \in 2^\omega : \exists^\infty n \ x \rest n+1 \in
A^n_{\stem(T)}\}$.
Clearly $\mu(B)\geq \delta_{n_0}$, so $B \cap X\neq \emptyset$.
Fix $x \in B\cap X$. We will find $S \geq T$ such that $S
\forces_{\Pg_g} x \in [\dot{A}]$, which will give a contradiction.

For each $n$ such that $x \in A^n_{\stem(T)}$ let $S_n \subseteq
T\rest n$ be a finite tree such that 
\begin{enumerate}
\item $\stem(S_n)=\stem(T)$,
\item for every $t \in S_n$, $n_0<|t|<n$, $\nor{\suc_{S_n}(t)} \geq
  \nor{\suc_T(t)}$,
\item for every $t \in S_n, \ |t|=n$, $x \in A_t$.
\end{enumerate}
The existence of $S_n$ follows from the inductive definition of
$A^n_t$'s.
By K\"onig lemma, there exists $S \subseteq T$ such that for
infinitely many $n$, $S \rest n=S_n$.
It follows that $S \in \Pg_g$ and
$S \forces_{\Pg_g} \exists^\infty n \ x\rest n \in \dot{A}\rest n$.
Since $ \dot{A}$ is a tree we conclude that $s \forces_{\Pg_g} x \in
[\dot{A}]$.~$\QED$

\section{More on $\UNIF({\mathcal J})$}
In this section we will discuss the model obtained by iterating the
forcing $\Pg_g$ alone.

\begin{theorem}
  It is consistent with $\ZFCa$ that for every $\sigma$-ideal
  ${\mathcal J}$ such that $\unif({\mathcal J})<2^{\boldsymbol\aleph_0}$, 
  $$\UNIF({\mathcal J} ) \subseteq \UNIF(\SN) \subseteq
  [\reals]^{<2^{\boldsymbol\aleph_0}}.$$
\end{theorem}
\Proof
Elements of $\UNIF(\SN)$ are traditionally called $C'$-sets. As we
remarked earlier, $\SN=\UNIF^\star(\SN)$. However, in
\cite{FreMil88Som} it is proved that assuming $\CH$, $\UNIF(\SN)
\subsetneq \SN $.

Let $\V$ be a model satisfying $\CH$ and let $\{g_\alpha:
\alpha<\omega_1\} \subseteq \omega^\omega $ be a dominating
family. Let $\{S_\alpha : \alpha < \omega_1\}$ be such that
\begin{enumerate}
\item $S_\alpha \cap S_\beta = \emptyset$ for $\alpha \neq \beta $,
\item $S_\alpha \subseteq \{\xi<\omega_2: \cf(\xi)=\omega_1\}$,
\item $S_\alpha $ is stationary for all $\alpha $.
\end{enumerate}

Let $\<{\mathcal P}_\alpha, \dot{{\mathcal
    Q}}_\alpha:\alpha<\omega_2\>$ be a countable support iteration
such that for $ \beta \in S_\alpha $, $\forces_\beta  \dot{{\mathcal
    Q}}_\beta  \simeq {\mathbb P}_{g_\alpha}$. If $ \beta \not \in
\bigcup_\alpha S_\alpha $ let $\dot{{\mathcal Q}}_\beta $ be trivial forcing.

Suppose that ${\mathcal J} $ is a $\sigma $-ideal and $\V^{{\mathcal
    P}_{\omega_2}} \thinks \unif({\mathcal J})=\boldsymbol\aleph_1 $. 
It follows that for some $\alpha<\omega_2$, $\V^{{\mathcal
    P}_{\omega_2}} \thinks \V^{{\mathcal P}_{\alpha}}\cap 2^\omega
\not \in {\mathcal J}$.

Suppose that $X \subseteq \V^{{\mathcal
    P}_{\omega_2}} \cap 2^\omega $ is uncountable. 

\bigskip 

{\sc Case 1} $|X|=\boldsymbol\aleph_1 $ and $\V^{{\mathcal
    P}_{\omega_2}} \thinks X \not\in \UNIF(\SN)$.

Let $\beta >\alpha$
 be such that 
\begin{enumerate}
\item $X \in \V^{{\mathcal
    P}_{\beta  }}$,
\item there is a continuous function $H: X \longrightarrow 2^\omega $, 
  $H \in \V^{{\mathcal P}_\beta }$ such that $\V^{{\mathcal P}_\beta}
    \thinks H"(X) \not\in \SN$,
\item $\beta \in S_\gamma $ and $\V^{{\mathcal
    P}_{\omega_2}} \thinks g_\gamma \text{ witnesses that } H"(X) \not\in
\UNIF(\SN)$.
\end{enumerate} 
It follows from the properties of $\Pg_{g_\gamma}$ that
$\V^{{\mathcal
    P}_{\beta+1}} \thinks \bigcup_n F_n"(H"(X)) =2^\omega  $. Hence
$\V^{{\mathcal
    P}_{\omega_2}} \thinks \bigcup_n F_n"(H"(X)) \not \in {\mathcal J} $,
which means that $\V^{{\mathcal
    P}_{\omega_2}} \thinks \exists n \in \omega  \ F_n"(H"(X)) \not \in
{\mathcal J}$.

\bigskip

{\sc Case 2} $|X|=2^{\boldsymbol\aleph_0} ={\mathbf\aleph}_2 $.

It is well known (see \cite{GJS92} or theorem 8.2.14 of \cite{BJbook})
in a model obtained by a countable support iteration of $\omega^\omega
$-bounding forcing notions there are no strong measure zero sets of
size $ 2^{\boldsymbol\aleph_0} $.
In particular, $\V^{{\mathcal P}_{\omega_2}} \thinks X \not \in
\SN$. Let $g_\beta $ be a witness to that.
Let $X_\beta  = X \cap \V^{{\mathcal P}_\beta}$. Standard argument
shows that
$$C=\{\gamma<\omega_2: \V^{{\mathcal P}_\gamma} \thinks X_\gamma
\not\in \SN\}$$
is a $\omega_1$-club.  Fix $\delta \in C \cap S_\beta $ and argue as
in the Case 1.~$\QED$

{\bf Acknowledgements}:
The work was done while the first author  was spending  his sabbatical year at 
 the 
Rutgers University and the College of Staten
Island, CUNY, and  their support is gratefully acknowledged.

\end{document}